\documentclass[final,1p,times]{elsarticle} 
\usepackage{amsmath}
\usepackage{hyphenat}
\newcommand{\dotminus}{\buildrel\textstyle.\over{\hbox{\vrule height2pt depth0pt width0pt}{\smash-}}}
\usepackage{amssymb}
\newtheorem{dummy}{Dummy}[section]
\newtheorem{definition}[dummy]{Definition}
\newtheorem{lemma}[dummy]{Lemma}
\newtheorem{corollary}[dummy]{Corollary}
\newtheorem{theorem}[dummy]{Theorem}

\journal{Annals of Pure and Applied Logic}

\begin{document}
\renewcommand\refname{}
\begin{frontmatter}

\title{Conditional computability of real functions \\with
respect to a class of operators}

\author[label1]{Ivan Georgiev}
\ead{ivandg@yahoo.com}
\author[label2]{Dimiter Skordev\corref{cor1}}
\ead{skordev@fmi.uni-sofia.bg}
\ead[url]{http://www.fmi.uni-sofia.bg/fmi/logic/skordev/}
\address[label1]{Burgas Prof.\@ Assen Zlatarov University, Faculty of Natural Sciences,\\1 Prof.\@ Yakimov blvd., 8010 Burgas, Bulgaria}
\address[label2]{Sofia University, Faculty of Mathematics and Informatics,\\5 James Bourchier blvd., 1164 Sofia, Bulgaria}
\cortext[cor1]{Corresponding author.}

\begin{abstract}
For any class of operators which transform unary total functions in the set of natural numbers into functions of the same kind, we define what it means for a real function to be uniformly computable or conditionally computable with respect to this class. These two computability notions are natural generalizations of certain notions introduced in a previous paper co-authored by Andreas Weiermann and in another previous paper by the same authors, respectively. Under certain weak assumptions about the class in question, we show that conditional computability is preserved by substitution, that all conditionally computable real functions are locally uniformly computable, and that the ones with compact domains are uniformly computable. The introduced notions have some similarity with the uniform computability and its non-uniform extension considered by Katrin Tent and Martin Ziegler, however, there are also essential differences between the conditional computability and the non-uniform computability in question.
\end{abstract}


\end{frontmatter}

\section{Introduction}\label{S0}
In the paper \cite{s-w-g}, a notion of uniform computability of real functions was introduced, namely, when a class $\cal F$ of total functions in $\mathbb{N}$ is given, some real functions were called uniformly $\cal F${\hyp}computable. The definition of the notion was in the spirit of the approach to computability of real functions originating from \cite{g,l} and nowadays indicated by the abbreviation TTE (cf., for instance, the monograph \cite{w}). This approach uses (necessarily infinitistic) naming systems for the real numbers and defines the computability of a real function as the existence of some effective procedure which transforms arbitrary names of the arguments into a name of the corresponding function value. The class of the computable real functions may depend on the choice of the naming system and on the sort of effective procedures which are admitted. However, there are some choices that hopefully produce the most general intuitively reasonable notion of computability for real functions (we will use the term {\em TTE computability} for this notion). Such a choice is, for instance, naming the real numbers by sequences of rational numbers converging to them with a given polynomial or exponential rate and transformation of the names through recursive operators (or, as for instance in \cite{w}, through oracle Turing machines). As indicated in~\cite{sh}, the restriction to general recursive operators leads to a narrower notion in the case of real functions which are not everywhere defined. Further reduction of the class of operators for the transformation of the names could additionally reduce the corresponding class of computable real functions, and this could be useful for introducing some subrecursive computability notions for real functions. 

In the case of uniform $\cal F${\hyp}computability, real numbers are named (up to technical details) by sequences of rational numbers converging to them with a linear rate, and the transformation of the names is performed by so-called $\cal F${\hyp}substitutional operators. Roughly speaking, the values of the image of a tuple of functions under such an operator are computed through evaluation of a term built from a variable, ranging over $\mathbb{N}$, by means of symbols for the functions in the tuple and for functions from $\cal F$. (In general, since $\cal F$ could contain some non{\hyp}computable functions, the corresponding procedure of transformation of names may be non-effective, and some uniformly $\cal F${\hyp}computable real functions may turn out to be not TTE computable. However, this cannot happen if $\cal F$ consists of recursive functions.)

The main attention in \cite{s-w-g} is paid to the case when $\cal F$ is a rather small subrecursive class, namely the class ${\cal M}^2$ (up to the argumentless constants $0,1,2,\ldots$, it consists of all functions in $\mathbb{N}$ which can be obtained from the successor function, the function $\lambda xy.x\dotminus  y$, the multiplication function and the projection functions by finitely many applications of substitution and bounded least number operation). Somewhat surprisingly, the results from \cite{s-w-g} easily imply that all elementary functions of calculus are uniformly ${\cal M}^2${\hyp}computable on the compact subsets of their domains. As to the uniform ${\cal M}^2${\hyp}computability of these functions on their whole domains, however, there is a serious obstacle for many of them, since any uniformly ${\cal M}^2${\hyp}computable real function is bounded by some polynomial.

In the paper \cite{s-g}, we introduced a wider notion of $\cal F${\hyp}computability called conditional $\cal F${\hyp}computability. Its definition ensures that all conditionally $\cal F${\hyp}computable real functions are TTE computable in the case when $\cal F$ consists of recursive functions. Under some weak assumptions on $\cal F$, we proved that conditional $\cal F${\hyp}computability is preserved by substitution, all conditionally $\cal F${\hyp}computable real functions are locally uniformly $\cal F${\hyp}computable and all conditionally $\cal F${\hyp}computable real functions with compact domains are uniformly $\cal F${\hyp}computable. Moreover, we prove that all elementary functions of calculus (considered on their whole domains) are conditionally ${\cal M}^2${\hyp}computable.\footnote{The proof in \cite{s-g} of the last statement needs some refinement concerning the functions $\lambda\xi.\sqrt[n]{\xi}$, $n=2,3,\ldots$ Although the exponential function and the logarithmic function are proved in that paper to be conditionally ${\cal M}^2${\hyp}computable, the expression for $\sqrt[n]{\xi}$ through them does not prove the conditional ${\cal M}^2${\hyp}computability of the functions $\lambda\xi.\sqrt[n]{\xi}$ on their whole domains (including 0 and also the negative real numbers in the case when $n$ is odd). The functions in question are actually uniformly ${\cal M}^2${\hyp}computable. This is proved for the case $n=2$ in \cite{s-w-g} and the proof can be easily modified to encompass the other values of $n$.} We also show the existence of TTE computable real functions which are not conditionally $\cal F$-computable, whatever be the class $\cal F$.

The supplementary feature of conditional $\cal F${\hyp}computability in comparison to the uniform one can be informally described as follows. It is now allowed the transformation of the names of the argument values (which produces a name for the corresponding function value) to depend on an additional parameter whose value is some natural number. Some term of the sort mentioned above must exist, such that whenever some names of the argument values are given, this number can be found by means of a search until the term in question vanishes (no restriction is imposed on the means used for organizing the search itself).\footnote{In the paper \cite{t-z}, two notions of computability of real functions with respect to some class~$\cal F$ of total functions in $\mathbb{N}$ are introduced and studied, namely the notion of a real function to be uniformly in $\cal F$ and the wider notion of a
real function to be in $\cal F$. The definitions of these notions use rational approximations of real numbers more directly than it is done in the TTE approach, but nevertheless there is a similarity to them of the uniform and the conditional computability, respectively. As seen from \cite{s}, the similarity of the uniform computability to the first notion is not a superficial one. However, the resemblance between conditional computability and the second notion is not so deep -- although the definition of this notion allows a dependence of the approximation process on an additional parameter, the description of its value uses the distance to the complement of the domain of the function, and the means provided by $\cal F$ could be not sufficient for checking if a given number is appropriate as a value of this parameter (a confusing feature of this notion is that some real functions with complicated domains may turn out to be in the class of the recursive functions without being TTE computable).} In order to take into account the parameter's value, we somewhat enlarged the class of procedures used for the transformation of names, and, roughly speaking, this corresponds to realizing them by using terms built in the above-mentioned way from two variables, ranging over $\mathbb{N}$, instead of one.

The present paper is devoted to a generalization of a part of the considerations in~\cite{s-g} to a situation when the class of the $\cal F${\hyp}substitutional operators is replaced with an arbitrary class $\mathbf{O}$ of operators in the set of the total unary functions in~$\mathbb{N}$. The dependence on the value of the additional parameter in the case of conditional computability is now realized by adding the corresponding constant function as an additional argument of the operators.

It seems that the approach described above has to do not only with the computational complexity of a real function, but also with the complexity of its definition. There is a more resource-oriented approach to computability of real functions due to Ker-I Ko in \cite{k}, where the complexity
of computable real functions is connected with discrete
polynomial complexity theory. 
The links between our
approach and Ko's one are yet to be specified.
At first glance, there is an essential difference between the two
approaches, namely Ko's approximations of the real numbers are defined to be of exponential rate,
but our approximations are necessarily of polynomial rate for the case, when
we compute by means of functions in $\mathbb{N}$, which are bounded by polynomial.
\vskip1cm

\section{Appropriate classes of operators and their relation to
$\cal F${\hyp}substitutional mappings}\label{S1}
\label{ac}
As we did in \cite{s-g}, for any $m\in\mathbb{N}$, we will denote by $\mathbb{T}_m$ the set of all
$m$-argument total functions in $\mathbb{N}$, and for any subset $\cal F$ of
$\bigcup_{m\in\mathbb{N}}\mathbb{T}_m$
and any
\mbox{$k,m\in\mathbb{N}$,} we may consider the notion of ${\cal F}${\hyp}substitutional mapping of
$\mathbb{T}_1^k$
into~$\mathbb{T}_m$. We will be interested in the last notion mainly for the case $m = 1$, and, in particular, we will indicate a way to reduce the case
$m > 1$ to this one.

In the present paper, the term ``operator'' will be used in the following restricted sense. For any $k\in\mathbb{N}$, the mappings of $\mathbb{T}_1^k$
into $\mathbb{T}_1$ (that is the $k$-ary operations in $\mathbb{T}_1$) will be called {\em $k$-ary operators}.
The $k$-ary operators for all $k\in\mathbb{N}$ will generally be called {\em operators}.

For any natural number $c$, let
$\check{c}$ be the constant function from $\mathbb{T}_1$ with
value~$c$. The identity function in $\mathbb{N}$ will be denoted by $\mathrm{id}_\mathbb{N}$. We set $\mathbb{T}=\bigcup_{m\in\mathbb{N}}\mathbb{T}_m$. We introduce the notion of appropriate class of operators, which, roughly speaking, captures
some necessary substitutional properties, possessed by the class of the $\cal F${\hyp}substitutional operators (as seen from Lemma \ref{L1}). 

\begin{definition}\label{D1}
Let $\mathbf{O}$ be a class of operators. The class $\mathbf{O}$ will be called
{\em appropriate}, whenever the following conditions are satisfied:
\begin{enumerate}
\item For any $k\in\mathbb{N}$ and any $i\in\{1,\ldots,k\}$, the operator $F$ defined by means of
the equality
$F(f_1,\ldots,f_k) = f_i$ 
belongs to $\mathbf{O}$.
\item The operator $F$ defined by
$F(f_1,f_2)(n) = f_1(f_2(n))$ 
belongs to $\mathbf{O}$.
\item For any $k,l\in\mathbb{N}$, if $F$ is a $k$-ary operator belonging to $\mathbf{O}$, and $G_1,\ldots,G_k$ are
$l$-ary operators belonging to $\mathbf{O}$, the operator $H$ defined by
$$H(g_1,\ldots,g_l) = F(G_1(g_1,\ldots,g_l),\ldots,G_k(g_1,\ldots,g_l))$$
also belongs to $\mathbf{O}$.
\item For any $k\in\mathbb{N}$ and any $(k+1)$-ary operator $F$ belonging to $\mathbf{O}$, the operator $G$
defined by
$$G(f_1,\ldots,f_k)(n) = F(f_1,\ldots,f_k,\check{n})(n)$$
also belongs to $\mathbf{O}$.
\end{enumerate}
\end{definition}

For the sake of convenience, we recall the definition from \cite{s-g} of the notion of $\cal F${\hyp}substitutional mapping.

\begin{definition}\label{substitutional}
Let ${\cal F}\subseteq\mathbb{T}$. For any $k,m\in\mathbb{N}$, certain mappings of $\mathbb{T}_1^k$ into $\mathbb{T}_m$ will be called {\em $\cal F${\hyp}substitutional}. We proceed by induction:
\begin{enumerate}
\item For any $m$-argument projection function $h$ in $\mathbb{N}$ the mapping $F$ defined by means of the equality $F(f_1,\ldots,f_k)=h$ is $\cal F${\hyp}substitutional.
\item For any $i\in\{1,\ldots,k\}$, if $F_0$ is an $\cal F${\hyp}substitutional mapping of $\mathbb{T}_1^k$ into $\mathbb{T}_m$ then so is the mapping $F$ defined by means of the equality
$$F(f_1,\ldots,f_k)(n_1,\ldots,n_m)=f_i(F_0(f_1,\ldots,f_k)(n_1,\ldots,n_m)).$$
\item For any natural number $r$ and any $r$-argument function $f$ from $\cal F$, if $F_1,\ldots,F_r$ are $\cal F${\hyp}substitutional mappings of $\mathbb{T}_1^k$ into $\mathbb{T}_m$ then so is the mapping $F$ defined by means of the equality
$$F(f_1,\ldots,f_k)(n_1,\ldots,n_m)=
f(F_1(f_1,\ldots,f_k)(n_1,\ldots,n_m),\ldots,F_r(f_1,\ldots,f_k)(n_1,\ldots,n_m)).$$ 
\end{enumerate}
\end{definition}

According to the terminology adopted in this paper, the $\cal F${\hyp}substitutional mappings of $\mathbb{T}_1^k$ into $\mathbb{T}_1$ will be called $k$-ary $\cal F${\hyp}substitutional operators.

\begin{lemma}\label{L1}
Let ${\cal F}\subseteq\mathbb{T}$
and let $\mathbf{O}$ be the class of all $\cal F${\hyp}substitutional
operators. Then $\mathbf{O}$ is appropriate.
\end{lemma}

{\bf Proof.} Condition 1 on $\mathbf{O}$ is satisfied by clauses 1 and 2 of Definition \ref{substitutional}
(considered for the case $m = 1$, when it coincides with the
corresponding Definition~6 from \cite{s-w-g}). Condition 2 on $\mathbf{O}$ follows from condition
1 on $\mathbf{O}$ and clause~2 of Definition \ref{substitutional}. Condition 3 on $\mathbf{O}$ is the case
$m = 1$ of Proposition~2 in \cite{s-g}. To prove condition 4 on $\mathbf{O}$ we use induction on
the construction of the operator $F$.

If $F$ is $\cal F${\hyp}substitutional by clause 1 of Definition \ref{substitutional}, then so is $G$.

Suppose that $F$ has the form from clause 2 of Definition \ref{substitutional}, that is
$$F(f_1,\ldots,f_k,f_{k+1})(n) = f_i(F_0(f_1,\ldots,f_k,f_{k+1})(n))$$ for some $i\in\{1,\ldots,k+1\}$ and, by
the inductive hypothesis, the operator $F_0$ has the considered property.
If $i\le k$, then
$$G(f_1,\ldots,f_k)(n) = F(f_1,\ldots,f_k,\check{n})(n) = f_i(F_0(f_1,\ldots,f_k,\check{n})(n)),$$
so by clause 2 of Definition \ref{substitutional}, $G$ is $\cal F${\hyp}substitutional.
If $i = k+1$, then
$$G(f_1,\ldots,f_k)(n) = F(f_1,\ldots,f_k,\check{n})(n) =
\check{n}(F_0(f_1,\ldots,f_k,\check{n})(n)) = n,$$
so by clause 1 of Definition \ref{substitutional}, $G$ is $\cal F${\hyp}substitutional.

Finally, suppose that $F$ is defined by
$$F(f_1,\ldots,f_k,f_{k+1})(n) = f(F_1(f_1,\ldots,f_k,f_{k+1})(n),\ldots,F_r(f_1,\ldots,f_k,f_{k+1})(n)),$$
where $f:\mathbb{N}^r\to\mathbb{N}$ belongs to $\cal F$ (clause 3 from Definition \ref{substitutional}) and by the
inductive hypothesis the operators $F_1,\ldots,F_r$ have the considered property. Then
$$G(f_1,\ldots,f_k)(n) = F(f_1,\ldots,f_k,\check{n})(n) =
f(F_1(f_1,\ldots,f_k,\check{n})(n),\ldots,F_r(f_1,\ldots,f_k,\check{n})(n)),$$ 
so $G$ is $\cal F${\hyp}substitutional by
clause 3 of Definition \ref{substitutional}.$\ \ \Box$

\begin{lemma}\label{L2}
Let ${\cal F}\subseteq\mathbb{T}$ and $F:\mathbb{T}_1^k\to\mathbb{T}_{m+1}$ for some $k,m\in\mathbb{N}$, $m\ne 0$.
The mapping $F$ is $\cal F${\hyp}substitutional iff there exists an $\cal F${\hyp}substitutional
mapping $G:\mathbb{T}_1^{k+1}\to\mathbb{T}_m$, such that
\begin{equation}\label{E1}
F(f_1,\ldots,f_k)(s,t_1,\ldots,t_m) = G(f_1,\ldots,f_k,\check{s})(t_1,\ldots,t_m)
\end{equation}
for all $f_1,\ldots,f_k\in\mathbb{T}_1$ and all $s,t_1,\ldots,t_m\in\mathbb{N}$.
\end{lemma}

{\bf Proof.} $(\Leftarrow)$. By induction on the construction of $G$, we will show that for any
$\cal F${\hyp}substitutional mapping $G:\mathbb{T}_1^{k+1}\to\mathbb{T}_m$, the mapping $F:\mathbb{T}_1^k\to\mathbb{T}_{m+1}$, defined
by the equality (\ref{E1}) is $\cal F${\hyp}substitutional.

If $G$ is $\cal F${\hyp}substitutional by clause 1 of Definition \ref{substitutional}, then so is $F$.

Let $G$ have the form
$$G(f_1,\ldots,f_k,f_{k+1})(t_1,\ldots,t_m) = f_i(G_0(f_1,\ldots,f_k, f_{k+1})(t_1,\ldots,t_m))$$
for some $i\in\{1,\ldots,k+1\}$ and let by the inductive hypothesis $G_0$ possess the
required property. If $i\le k$, then
$$F(f_1,\ldots,f_k)(s,t_1,\ldots,t_m) =
f_i(G_0(f_1,\ldots,f_k,\check{s})(t_1,\ldots,t_m)),$$
so $F$ is $\cal F${\hyp}substitutional
by clause 2 of Definition \ref{substitutional}. If $i = k+1$, then
$$F(f_1,\ldots,f_k)(s,t_1,\ldots,t_m) =
\check{s}(G_0(f_1,\ldots,f_k,\check{s})(t_1,\ldots,t_m)) = s,$$
so $F$ is $\cal F${\hyp}substitutional by clause 1 of Definition \ref{substitutional}.

Finally, suppose that $G$ is defined by
\begin{multline*}
G(f_1,\ldots,f_k,f_{k+1})(t_1,\ldots,t_m)\\ =
f(G_1(f_1,\ldots,f_k,f_{k+1})(t_1,\ldots,t_m),\ldots,G_r(f_1,\ldots,f_k, f_{k+1})(t_1,\ldots,t_m))
\end{multline*}
for a function $f:\mathbb{N}^r\to\mathbb{N}$ belonging to $\cal F$ (clause 3 of Definition \ref{substitutional}) and
mappings $G_1,\ldots,G_r$, which by the inductive hypothesis possess the
considered property. Then
$$F(f_1,\ldots,f_k)(s, t_1,\ldots,t_m) = 
f(G_1(f_1,\ldots,f_k,\check{s})(t_1,\ldots,t_m),\ldots,G_r(f_1,\ldots,f_k,\check{s})(t_1,\ldots,t_m)),$$ 
so $F$ is
$\cal F${\hyp}substitutional by clause 3 of Definition \ref{substitutional}.

$(\Rightarrow)$. By induction on the construction of $F$, we will show that for any
$\cal F${\hyp}substitutional mapping $F:\mathbb{T}_1^k\to\mathbb{T}_{m+1}$, there exists an $\cal F${\hyp}substitutional
mapping $G:\mathbb{T}_1^{k+1}\to\mathbb{T}_m$, such that the equality (\ref{E1}) holds for all
$f_1,\ldots,f_k\in\mathbb{T}_1$ and all $s,t_1,\ldots,t_m\in\mathbb{N}$.

Suppose $F(f_1,\ldots,f_k)(s,t_1,\ldots,t_m) = t_i$ for some $i\in\{ 1,\ldots,m\}$.
Then the mapping $G$, defined by $G(f_1,\ldots,f_k,f_{k+1})(t_1,\ldots,t_m) = t_i$ satisfies the
required condition.

If $F(f_1,\ldots,f_k)(s,t_1,\ldots,t_m) = s$, then the mapping $G$, defined by
$$G(f_1,\ldots,f_k,f_{k+1})(t_1,\ldots,t_m) = f_{k+1}(t_1)$$ is ${\cal F}${\hyp}substitutional (by clauses 1 and 2 of
Definition \ref{substitutional}) and
$$G(f_1,\ldots,f_k,\check{s})(t_1,\ldots,t_m) = \check{s}(t_1) = s = F(f_1,\ldots,f_k)(s,t_1,\ldots,t_m).$$

If $F$ has the form
$$F(f_1,\ldots,f_k)(s,t_1,\ldots,t_m) = f_i(F_0(f_1,\ldots,f_k)(s,t_1,\ldots,t_m))$$
for some $i\in\{1,\ldots,k\}$ and $G_0$ is the $\cal F${\hyp}substitutional mapping, which exists
by the inductive hypothesis for the mapping $F_0$, then we can define a
mapping $G$ by $$G(f_1,\ldots,f_k,f_{k+1})(t_1,\ldots,t_m) = f_i(G_0(f_1,\ldots,f_k,f_{k+1})(t_1,\ldots,t_m)),$$ 
which is $\cal F${\hyp}substitutional by clause 2 of Definition \ref{substitutional}. It follows that
$$F(f_1,\ldots,f_k)(s, t_1,\ldots,t_m) = 
f_i(G_0(f_1,\ldots,f_k,\check{s})(t_1,\ldots,t_m)) = G(f_1,\ldots,f_k,\check{s})(t_1,\ldots,t_m).$$ 

Finally, suppose that $F$ is defined by
$$F(f_1,\ldots,f_k)(s,t_1,\ldots,t_m) = 
f(F_1(f_1,\ldots,f_k)(s,t_1,\ldots,t_m),\ldots,F_r(f_1,\ldots,f_k)(s, t_1,\ldots,t_m)),$$ 
where $f:\mathbb{N}^r\to\mathbb{N}$ belongs to $\cal F$ (clause 3 from Definition \ref{substitutional}) and by the
inductive hypothesis the mappings $F_1,\ldots,F_r$ have the considered property,
that is there exist $\cal F${\hyp}substitutional mappings $G_1,\ldots,G_r$, such that
$$F_i(f_1,\ldots,f_k)(s,t_1,\ldots,t_m) = G_i(f_1,\ldots,f_k,\check{s})(t_1,\ldots,t_m)$$ for all $i\in\{1,\ldots,r\}$.
Let us define
\begin{multline*}
G(f_1,\ldots,f_k,f_{k+1})(t_1,\ldots,t_m)\\ =
f(G_1(f_1,\ldots,f_k,f_{k+1})(t_1,\ldots,t_m),\ldots,G_r(f_1,\ldots,f_k, f_{k+1})(t_1,\ldots,t_m)).
\end{multline*}
Then $G$ is $\cal F${\hyp}substitutional by clause 3 of Definition \ref{substitutional} and
\begin{multline*}
F(f_1,\ldots,f_k)(s, t_1,\ldots,t_m) =
f(G_1(f_1,\ldots,f_k,\check{s})(t_1,\ldots,t_m),\ldots,G_r(f_1,\ldots,f_k,\check{s})(t_1,\ldots,t_m))\\ =
G(f_1,\ldots,f_k,\check{s})(t_1,\ldots,t_m).\ \ \Box
\end{multline*}

\begin{corollary}\label{C1}
Let ${\cal F}\subseteq\mathbb{T}$
and $F:\mathbb{T}_1^k\to\mathbb{T}_{m+1}$ for some $k,m\in\mathbb{N}$.
The mapping $F$ is $\cal F${\hyp}substitutional iff there exists an $\cal F${\hyp}substitutional
operator $G:\mathbb{T}_1^{k+m}\to\mathbb{T}_1$, such that
$$F(f_1,\ldots,f_k)(s_1,\ldots,s_m,t) = G(f_1,\ldots,f_k,\check{s_1},\ldots,\check{s_m})(t)$$
for all $f_1,\ldots,f_k\in\mathbb{T}_1$ and all $s_1,\ldots,s_m,t\in\mathbb{N}$.
\end{corollary}

For any function $f:\mathbb{N}^k\to\mathbb{N}$ we define the $k$-ary operator
$\mathring{f}$ by the equality
$$\mathring{f}(f_1,\ldots,f_k)(n) = f(f_1(n),\ldots,f_k(n))$$
for all $n\in\mathbb{N}$.

\begin{definition}\label{D2}
Let $\mathbf{O}$ be a class of operators. A function $f:\mathbb{N}^k\to\mathbb{N}$ will be
called {\em representable} in $\mathbf{O}$ or $\mathbf{O}$-{\em representable}, for short, if the
corresponding operator
$\mathring{f}$ belongs to $\mathbf{O}$.
\end{definition}

\begin{lemma}\label{L3}
Let ${\cal F}\subseteq\mathbb{T}$
and $\mathbf{O}$ be the class of all $\cal F${\hyp}substitutional
operators. Then all functions from $\cal F$ are $\mathbf{O}$-representable.
\end{lemma}

{\bf Proof.} We use clause 3 in Definition \ref{substitutional} and condition 1 in Definition \ref{D1},
making use of the fact that $\mathbf{O}$ is appropriate.$\ \ \Box$

\begin{lemma}\label{L4}
Let $\mathbf{O}$ be an appropriate class of operators.
Then all projection functions in $\mathbb{N}$ are $\mathbf{O}$-representable.
\end{lemma}

{\bf Proof.} Let $f\!:\!\mathbb{N}^k\!\!\to\!\!\mathbb{N}$ be a projection function, that is for some $i\!\in\!\{1,\ldots,k\}$, 
\mbox{$f(n_1,\ldots,n_k) = n_i$.} Then 
$\mathring{f}(f_1,\ldots,f_k)(n) = f(f_1(n),\ldots,f_k(n)) = f_i(n)$ for all $n\in\mathbb{N}$,
that is
$\mathring{f}(f_1,\ldots,f_k) = f_i$ and by condition 1 in Definition \ref{D1},
$\mathring{f}\in\mathbf{O}$.
It follows that $f$ is $\mathbf{O}$-representable.$\ \ \Box$

\begin{lemma}\label{L5}
Let $\mathbf{O}$ be an appropriate class of operators. Then the class of all
$\mathbf{O}$-representable functions is closed under substitution.
\end{lemma}

{\bf Proof.} Let $f:\mathbb{N}^k\to\mathbb{N}$ and $g_i:\mathbb{N}^l\to\mathbb{N}$, $i = 1,\ldots,k$ be $\mathbf{O}$-representable functions
and the function $g:\mathbb{N}^l\to\mathbb{N}$ be defined by
$$g(n_1,\ldots,n_l) = f(g_1(n_1,\ldots,n_l),\ldots,g_k(n_1,\ldots,n_l)).$$
By condition 3 in Definition \ref{D1}, it is sufficient to show that
$$\mathring{g}(f_1,\ldots,f_l) =
\mathring{f}(\mathring{g_1}(f_1,\ldots,f_l),\ldots,
\mathring{g_k}(f_1,\ldots,f_l))$$ holds for all $f_1,\ldots,f_l\in\mathbb{T}_1$.
This can be seen as follows: for all $n\in\mathbb{N}$ and $f_1,\ldots,f_l\in\mathbb{T}_1$,
\begin{multline*}
\mathring{g}(f_1,\ldots,f_l)(n) = g(f_1(n),\ldots,f_l(n)) = f(g_1(f_1(n),\ldots,f_l(n)),\ldots,g_k(f_1(n),\ldots,f_l(n)))\\ 
= f(\mathring{g_1}(f_1,\ldots,f_l)(n),\ldots,\mathring{g_k}(f_1,\ldots,f_l)(n)) = \mathring{f}(\mathring{g_1}(f_1,\ldots,f_l),\ldots,\mathring{g_k}(f_1,\ldots,f_l))(n).\ \ \Box 
\end{multline*}

\vskip-3mm
\begin{lemma}\label{is}
Let $\mathbf{O}$ be an appropriate class of operators, and $k$ be a natural number. Then:
\begin{enumerate}
\item The $k$-ary operator $F$ defined by
$F(f_1,\ldots,f_k)=\mathrm{id}_\mathbb{N}$ 
belongs to $\mathbf{O}$.
\item For any $i\in\{1,\ldots,k\}$, if $F_0$ is a $k$-ary operator belonging to $\mathbf{O}$ then so is the operator $F$ defined by means of the equality
$F(f_1,\ldots,f_k)(n)=f_i(F_0(f_1,\ldots,f_k)(n))$. 
\item For any natural number $r$ and any $r$-argument $\mathbf{O}$-representable function~$f$, if $F_1,\ldots,F_r$ are $k$-ary operators belonging to $\mathbf{O}$ then so is the operator~$F$ defined by means of the equality
$F(f_1,\ldots,f_k)(n)=f(F_1(f_1,\ldots,f_k)(n),\ldots,F_r(f_1,\ldots,f_k)(n))$. 
\end{enumerate}
\end{lemma}

{\bf Proof.} If $F$ is the operator from statement 1 of the lemma then, for all $f_1,\ldots,f_k\in\mathbb{T}_1$ and all $n\in\mathbb{N}$,
$$F(f_1,\ldots,f_k)(n)=n=F_0(f_1,\ldots,f_k,\check{n})(n),$$
where the $(k+1)$-ary operator $F_0$ is defined by
$F_0(f_1,\ldots,f_k,f_{k+1})=f_{k+1}$, 
hence $F\in\mathbf{O}$ by conditions 1 and 4 of Definition \ref{D1}.

If $F$ is defined in the way from statement 2 then, for all $f_1,\ldots,f_k\in\mathbb{T}_1$,
$$F(f_1,\ldots,f_k)=F_1(f_i,F_0(f_1,\ldots,f_k)),$$
where $F_1$ is the operator from condition 2 of Definition \ref{D1}, hence $F\in\mathbf{O}$ by conditions 1, 2 and 3 of Definition \ref{D1}.

Finally, suppose that $F$ is defined in the way from statement 3. Then, for all $f_1,\ldots,f_k\in\mathbb{T}_1$,
$$F(f_1,\ldots,f_k)=\mathring{f}(F_1(f_1,\ldots,f_k),\ldots,F_r(f_1,\ldots,f_k)),$$
hence $F\in\mathbf{O}$ by condition 3 of Definition \ref{D1} and the fact that $f$ is $\mathbf{O}${\hyp}representable.$\ \ \Box$ 

\begin{corollary}\label{L7}
Let $\mathbf{O}$ be an appropriate class of operators, $k$ be a natural number, and $f$ be a unary $\mathbf{O}$-representable function. Then the operator $F$ defined by
$F(f_1,\ldots,f_k)=f$ 
belongs to $\mathbf{O}$.
\end{corollary}

\vskip3mm
{\bf Proof.} Let $F_1$ be the operator defined by
$F_1(f_1,\ldots,f_k)=\mathrm{id}_\mathbb{N}$. 
Then, for all $f_1,\ldots,f_k\in\mathbb{T}_1$ and all $n\in\mathbb{N}$,
$F(f_1,\ldots,f_k)(n)=f(F_1(f_1,\ldots,f_k)(n))$, 
hence $F\in\mathbf{O}$ by statements 1 and 3 of Lemma \ref{is}.$\ \ \Box$

\vskip3mm
\begin{corollary}\label{swrtOr}
Let $\mathbf{O}$ be an appropriate class of operators, and $\cal F$ be the class of all $\mathbf{O}$-representable functions. Then all $\cal F${\hyp}substitutional operators belong to $\mathbf{O}$.
\end{corollary}

{\bf Proof.} By induction making use of statements 1, 2 and 3 of Lemma \ref{is}.$\ \ \Box$

\vskip3mm
The continuity notion for operators is defined in the usual way.
\begin{definition}\label{D3}
Let $F:\mathbb{T}_1^k\to\mathbb{T}_1$. The operator $F$ will be called {\em continuous},
if the following condition is satisfied: for all $f_1,\ldots,f_k\in\mathbb{T}_1$ and $n\in\mathbb{N}$,
there exists $u\in\mathbb{N}$, such that
\begin{equation}\label{cont} 
F(g_1,\ldots,g_k)(n) = F(f_1,\ldots,f_k)(n),
\end{equation} 
whenever
$g_1,\ldots,g_k\in\mathbb{T}_1$ and $g_1(t) = f_1(t),\,\ldots,\,g_k(t) = f_k(t)$ for all $t\le u$.
\end{definition}

\begin{lemma}\label{L8}
For any ${\cal F}\subseteq\mathbb{T}$,
all $\cal F${\hyp}substitutional operators are
continuous.
\end{lemma}

{\bf Proof.} One shows by induction on the construction of $F$ that
every $\cal F${\hyp}substitutional operator $F$ is continuous.$\ \ \Box$

\vskip2mm
Clearly, the class of all operators is an appropriate one. By the above lemma, this class is different from the class of the $\cal F${\hyp}substitutional operators for any choice of the class ${\cal F}\subseteq\mathbb{T}$. There are, however, more interesting examples of appropriate classes with this property. Such classes are, for instance, the class of all computable operators, the class of the primitive recursive ones, the class of the elementary ones, etc. Although  consisting of continuous operators, each of these classes contains some operators which are not $\cal F${\hyp}substitutional, whatever be the class $\cal F$. This can be seen by inductively proving the stronger continuity property of the $\cal F${\hyp}substitutional operators formulated in the lemma below.

\begin{lemma}\label{stronger continuity}
Let $\cal F$ be an arbitrary subclass of $\mathbb{T}$. Then, whenever $F$ is a $k$-ary $\cal F${\hyp}substitutional operator, there exists a natural number $v$ with the following property: for all $f_1,\ldots,f_k\in\mathbb{T}_1$ and $n\in\mathbb{N}$,
there exists a subset $A$ of~$\mathbb{N}$, such that $A$ has at most $v$ elements, and the equality (\ref{cont}) holds, whenever 
$g_1,\ldots,g_k\in\mathbb{T}_1$ and $g_1(t) = f_1(t),\,\ldots,\,g_k(t) = f_k(t)$ for all $t\in A$.
\end{lemma}

We emphasize that the property from Lemma \ref{stronger continuity} is indeed stronger than ordinary continuity,
since the set $A$ may depend on the choice of $f_1,\ldots,f_k\in\mathbb{T}_1$ and $n\in\mathbb{N}$, but the cardinality bound~$v$ for $A$ may not. 

\section{Uniform and conditional computability of a real function \\with respect to a class of operators}\label{S2}

As in \cite{s-w-g}, a triple $(f,g,h)\in\mathbb{T}_1^3$
is called to {\em name} a real number $\xi$ if
$$\left|\frac{f(t)-g(t)}{h(t)+1}-\xi\right|<\frac{1}{t+1}$$
for all $t\in\mathbb{N}$.

\begin{definition}\label{D4}
Let $\mathbf{O}$ be a class of operators, $N\in\mathbb{N}$ and $\theta:D\to\mathbb{R}$, where $D\subseteq\mathbb{R}^N$. The function $\theta$ will be called {\em uniformly computable} with respect to~$\mathbf{O}$
or {\em uniformly $\mathbf{O}${\hyp}computable}, for short, if there exist $3N$-ary operators $F,G,H$
belonging to $\mathbf{O}$, such that whenever $(\xi_1,\ldots,\xi_N)\in D$ and $(f_1,g_1,h_1),\,\ldots,\,(f_N,g_N,h_N)$ are triples from $\mathbb{T}_1^3$ 
naming $\xi_1,\,\ldots,\,\xi_N$, respectively,
the triple
$$(F(f_1,g_1,h_1,\ldots,f_N,g_N,h_N),G(f_1,g_1,h_1,\ldots,f_N,g_N,h_N),H(f_1,g_1,h_1,\ldots,f_N,g_N,h_N))$$ 
names $\theta(\xi_1,\ldots,\xi_N)$.
\end{definition}

If ${\cal F}\subseteq\mathbb{T}$
and $\mathbf{O}$ is the class of all $\cal F${\hyp}substitutional operators,
then a real function is uniformly $\mathbf{O}${\hyp}computable iff it is
uniformly $\cal F${\hyp}computable in the sense of Definition 7 in \cite{s-w-g}.

\begin{definition}\label{D5}
Let $\mathbf{O}$ be a class of operators, $N\in\mathbb{N}$ and $\theta:D\to\mathbb{R}$, where \mbox{$D\subseteq\mathbb{R}^N$.} The function $\theta$ will be called {\em conditionally computable}
with respect to~$\mathbf{O}$ or {\em conditionally $\mathbf{O}${\hyp}computable}, for short, if there exist
a $3N$-ary operator~$E$ and $(3N+1)$-ary operators $F,G,H$, such that
$E,F,G,H\in\mathbf{O}$ and, whenever $(\xi_1,\ldots,\xi_N)\in D$ and $(f_1,g_1,h_1),\,\ldots,\,(f_N,g_N,h_N)$ are triples from $\mathbb{T}_1^3$ naming $\xi_1,\,\ldots,\,\xi_N$, respectively,
the following holds:
\begin{enumerate}
\item There exists a natural number $s$ satisfying the equality
\begin{equation}\label{E2}
E(f_1,g_1,h_1,\ldots,f_N,g_N,h_N)(s) = 0.
\end{equation}
\item For any natural number s satisfying the equality (\ref{E2}), the triple
\begin{multline*}
(F(f_1,g_1,h_1,\ldots,f_N,g_N,h_N,\check{s}),G(f_1,g_1,h_1,\ldots,f_N,g_N,h_N,\check{s}),\\H(f_1,g_1,h_1,\ldots,f_N,g_N,h_N,\check{s}))
\end{multline*}
names $\theta(\xi_1,\ldots,\xi_N)$.
\end{enumerate}
\end{definition}

By the case $m = 1$ of Lemma \ref{L2}, if ${\cal F}\subseteq\mathbb{T}$,
and $\mathbf{O}$ is the class of all
$\cal F${\hyp}substitutional operators, then a real function is conditionally
$\mathbf{O}${\hyp}computable iff it is conditionally $\cal F${\hyp}computable in the sense of
Definition 2 in \cite{s-g}.

If $\mathbf{O}$ is an appropriate class of operators
then all uniformly $\mathbf{O}${\hyp}computable real functions are conditionally $\mathbf{O}${\hyp}computable. Indeed,
let $\mathbf{O}$ be an appropriate class of operators and $F^\circ,G^\circ,H^\circ$ be the operators from Definition \ref{D4} belonging to $\mathbf{O}$ for the
$N$-argument real function $\theta$. Then we can satisfy the requirements from
Definition \ref{D5} by setting
\begin{align*}
E(f_1,g_1,h_1,\ldots,f_N,g_N,h_N)= &\,\mathrm{id}_\mathbb{N},\\
F(f_1,g_1,h_1,\ldots,f_N,g_N,h_N,e)\! = &\,F^\circ(f_1,g_1,h_1,\ldots,f_N,g_N,h_N),\\
G(f_1,g_1,h_1,\ldots,f_N,g_N,h_N,e)\! = &\,G^\circ(f_1,g_1,h_1,\ldots,f_N,g_N,h_N),\\
H(f_1,g_1,h_1,\ldots,f_N,g_N,h_N,e)\! = &\,H^\circ(f_1,g_1,h_1,\ldots,f_N,g_N,h_N).
\end{align*}
The operator $E$ belongs to $\mathbf{O}$ by statement 1 of Lemma \ref{is}. By conditions 1 and~3
from Definition \ref{D1} and by the fact that $F^\circ,G^\circ,H^\circ$ belong to $\mathbf{O}$, the operators
$F,G,H$ also belong to $\mathbf{O}$.

If $\mathbf{O}$ is a class of recursive operators then obviously all uniformly $\mathbf{O}${\hyp}computable real functions are TTE computable, and it is easy to show the same for the conditionally $\mathbf{O}${\hyp}computable ones.

\section{Substitution in conditionally $\mathbf{O}${\hyp}computable real functions}\label{S3}

The next theorem generalizes Theorem 1 of \cite{s-g}.

\begin{theorem}\label{T1}
Let $\mathbf{O}$ be an appropriate class of operators. Let a two-argument $\mathbf{O}$-representable function $C$ and $\mathbf{O}$-representable one-argument
functions $L$ and $R$ in $\mathbb{N}$ exist such that
$$\{(u,v)\in\mathbb{N}^2\,|\,C(u,v)=0\}=\{(0,0)\},\ \ \{(L(s),R(s))\,|\,s\in\mathbb{N}\} = \mathbb{N}^2.$$
Then the substitution operation on real functions preserves conditional
$\mathbf{O}${\hyp}computability.
\end{theorem}

{\bf Proof.} To avoid writing excessively long expressions, we will restrict ourselves
to the case of one-argument functions. Let $\theta_0$ and $\theta_1$ be conditionally
$\mathbf{O}${\hyp}computable one-argument real functions. We will show the conditional
$\mathbf{O}${\hyp}computability of the function $\theta$, defined by $$\theta(\xi) = \theta_0(\theta_1(\xi)).$$ For $i = 0,1$, 
let $E_i,F_i,G_i,H_i$ be operators from $\mathbf{O}$, such that $\exists s(E_i(f,g,h)(s) = 0)$ and
$$\forall s(E_i(f,g,h)(s) = 0 
\Rightarrow(F_i(f,g,h,\check{s}),G_i(f,g,h,\check{s}),H_i(f,g,h,\check{s}))\
\mathrm{names}\ \theta_i(\xi))$$ for any $\xi\in\mathrm{dom}(\theta_i)$ and any triple $(f,g,h)$ naming $\xi$. We will show that the
requirements of Definition \ref{D5} for the function $\theta$ are satisfied through the
operators $E,F,G,H$ defined as follows:\footnote{The function $C$ will be used to model conjunction. The functions $L$ and $R$ will be used for decoding a natural number~$s$, which codes the values of the parameters $s_0$ and $s_1$, corresponding to the two real functions $\theta_0$ and $\theta_1$, respectively.} 
\begin{equation*}\begin{split}
E(f,g,h)(s)\!=&C(E_1(f,g,h)(R(s)),E_0(F_1(f,g,h,\mathring{R}(\check{s})),G_1(f,g,h,\mathring{R}(\check{s})),H_1(f,g,h,\mathring{R}(\check{s})))(L(s))),\\ F(f,g,h,e)\!= &F_0(F_1(f,g,h,\mathring{R}(e)),G_1(f,g,h,\mathring{R}(e)),H_1(f,g,h,\mathring{R}(e)),\mathring{L}(e)),\\
G(f,g,h,e)\!= &G_0(F_1(f,g,h,\mathring{R}(e)),G_1(f,g,h,\mathring{R}(e)),H_1(f,g,h,\mathring{R}(e)),\mathring{L}(e)),\\
H(f,g,h,e)\!= &H_0(F_1(f,g,h,\mathring{R}(e)),G_1(f,g,h,\mathring{R}(e)),H_1(f,g,h,\mathring{R}(e)),\mathring{L}(e)).
\end{split}\end{equation*}
Suppose $\xi\in\mathrm{dom}(\theta)$ and $(f,g,h)$ is a triple naming $\xi$. By the conditional
$\mathbf{O}${\hyp}computability of $\theta_1$, there exists $s_1\in\mathbb{N}$ such that
\begin{equation}\label{E3}
E_1(f,g,h)(s_1) = 0,
\end{equation}
and if we choose such an $s_1$, then the number $\theta_1(\xi)$ is named
by the triple $(f_1,g_1,h_1)$, where
\begin{equation}\label{E4}
f_1 = F_1(f,g,h,\check{s_1}),\ g_1 = G_1(f,g,h,\check{s_1}),\ h_1 = H_1(f,g,h,\check{s_1}).
\end{equation}
By the conditional $\mathbf{O}${\hyp}computability of $\theta_0$, there exists $s_0\in\mathbb{N}$, such that
\begin{equation}\label{E5}
E_0(f_1,g_1,h_1)(s_0) = 0.
\end{equation}
If $s$ is a natural number, such that $L(s) = s_0$, $R(s) = s_1$, then $E(f,g,h)(s) = 0$.
Consider now any natural number $s$, such that $E(f,g,h)(s) = 0$.
Let $s_0 = L(s)$ and $s_1 = R(s)$. The equality $E(f,g,h)(s) = 0$ implies the
equality (\ref{E3}), as well as the equality (\ref{E5}) for the functions $f_1,g_1,h_1$,
defined by means of the equalities~(\ref{E4}). It follows from the equality (\ref{E3}) that
$(f_1,g_1,h_1)$ names $\theta_1(\xi)$, and, together with the equality (\ref{E5}), this fact implies
that $\theta(\xi) = \theta_0(\theta_1(\xi))$ is named by the triple
$$(F_0(f_1,g_1,h_1,\mathring{L}(\check{s})),G_0(f_1,g_1,h_1,\mathring{L}(\check{s})),H_0(f_1,g_1,h_1,\mathring{L}(\check{s}))),$$
which coincides with the triple
$$(F(f,g,h,\check{s}),G(f,g,h,\check{s}),H(f,g,h,\check{s})).$$
The operators $F,G,H$ belong to $\mathbf{O}$ by conditions 1 and 3 of Definition \ref{D1} and
the fact that all the operators $F_0,G_0,H_0,F_1,G_1,H_1,\mathring{R},\mathring{L}$
belong to $\mathbf{O}$.
It remains to show that the operator $E$ also belongs to $\mathbf{O}$.
Let the operators $A$ and $B$ be defined by
\begin{gather*}
A(f,g,h)(s) = E_1(f,g,h)(R(s)),\\
B(f,g,h)(s) = E_0(F_1(f,g,h,\mathring{R}(\check{s})),G_1(f,g,h,\mathring{R}(\check{s})),H_1(f,g,h,\mathring{R}(\check{s})))(L(s)).
\end{gather*}
Then we have
$$E(f,g,h)(s)=C(A(f,g,h)(s),B(f,g,h)(s)),$$
hence, by the fact that the function $C$ is $\mathbf{O}$-representable and by applying
statement 3 from Lemma \ref{is}, it would be sufficient to show that $A,B\in\mathbf{O}$
in order to conclude that $E\in\mathbf{O}$. Let $U:\mathbb{T}_1^3\to\mathbb{T}_1$ be defined by
$$U(f,g,h) = R$$ for all $f,g,h\in\mathbb{T}_1$
and $V:\mathbb{T}_1^4\to\mathbb{T}_1$
be defined by
$$V(f,g,h,e) = L$$ for all $f,g,h,e\in\mathbb{T}_1$.
The operators $U$ and $V$ belong to $\mathbf{O}$ thanks to Corollary~\ref{L7} and the fact that $R$ and $L$ are $\mathbf{O}$-representable. Then we have
$$A(f,g,h)(s) = E_1(f,g,h)(U(f,g,h)(s)),$$
so $A\in\mathbf{O}$ by conditions 2 and 3 from Definition \ref{D1} and the fact that $E_1,U\in\mathbf{O}$.
To show that $B$ also belongs to $\mathbf{O}$, we note that, for all $f,g,h\in\mathbb{T}_1$ and $s\in\mathbb{N}$,
$$B(f,g,h)(s) = W(f,g,h,\check{s})(s),$$
where the operator $W$ is defined by
$$W(f,g,h,e)(s) = E_0(F_1(f,g,h,\mathring{R}(e)),G_1(f,g,h,\mathring{R}(e)),H_1(f,g,h,\mathring{R}(e)))(L(s)).$$
Since $L(s)=V(f,g,h,e)(s)$, the operator $W$ belongs to $\mathbf{O}$ by conditions 1, 2 and 3 from Definition~\ref{D1} and 
the fact that $E_0,F_1,G_1,H_1,\mathring{R},V$ belong to $\mathbf{O}$. Therefore, by condition 4 from
Definition~\ref{D1}, the operator $B$ also belongs to~$\mathbf{O}$.$\ \ \Box$ 

\section{Local uniform $\mathbf{O}${\hyp}computability of the conditionally $\mathbf{O}${\hyp}computable real
functions}\label{S4}

\begin{definition}\label{D6}
Let $\mathbf{O}$ be a class of operators, $N\in\mathbb{N}$ and $\theta:D\to\mathbb{R}$,
where $D\subseteq\mathbb{R}^N$. The function $\theta$ will be called {\em locally uniformly $\mathbf{O}${\hyp}computable},
if any point of $D$ has some neighbourhood~$U$, such that the restriction of $\theta$ 
to $D\cap U$ is uniformly $\mathbf{O}${\hyp}computable.
\end{definition}

For any $k,c\in\mathbb{N}$, let the function $\mu_{k,c}:\mathbb{N}^2\to\mathbb{N}$ be defined as follows:\footnote{We will use the functions $\mu_{k,c}$ to model a particular nested if-then-else construction.}
$$\mu_{k,c}(x,y)=\left\{\begin{array}{ll}\!c & \mathrm{if}\ x=k,\\\!y & \mathrm{otherwise.}\end{array}\right.$$

The next theorem generalizes Theorem 2 of \cite{s-g}.

\begin{theorem}\label{T2}
Let $\mathbf{O}$ be an appropriate class of continuous operators, and let the functions $\check{c}$ for all $c\in\mathbb{N}$ and the functions $\mu_{k,c}$ for all $k,c\in\mathbb{N}$ be $\mathbf{O}$-representable. Then all conditionally $\mathbf{O}${\hyp}computable real functions are locally
uniformly $\mathbf{O}${\hyp}computable.
\end{theorem}

{\bf Proof.} For an arbitrary function $a\in\mathbb{T}_1$ and any $k\in\mathbb{N}$, let the unary operator ${^ka}$ be defined as follows:
$${^ka}(f)(t)=\left\{\begin{array}{ll}\!a(t) & \mathrm{if}\ t<k,\\\!f(t) & \mathrm{otherwise.}\end{array}\right.$$
We will show by induction on $k$ that $^ka\in\mathbf{O}$. The operator $^0a$ belongs to~$\mathbf{O}$ by clause 1 of Definition \ref{D1} since ${^0a}(f)=f$ for all $f\in\mathbb{T}_1$. Suppose now, by the inductive hypothesis, that $^ka\in\mathbf{O}$ for a certain $k\in\mathbb{N}$. Then $^{k+1}a$ also belongs to $\mathbf{O}$ by the equality
$${^{k+1}a}(f)(t)=\mu_{k,a(k)}\!\left(t,{^ka}(f)(t)\right),$$
statements 1 and 3 of Lemma \ref{is}, the inductive hypothesis and the $\mathbf{O}${\hyp}representability of $\mu_{k,a(k)}$.

Let now $\theta:D\to\mathbb{R}$, where $D\subseteq R$, be a conditionally $\mathbf{O}${\hyp}computable real function,
and $\xi_0\in D$ (for the sake of simplicity, we assume additionally that $\theta$ is
unary). Let $E,F,G,H\in\mathbf{O}$ be witnesses from Definition \ref{D5} (with $N = 1$).
Let $(f_0,g_0,h_0)$ be a triple naming $\xi_0$, and let $s_0$ be a natural number,
satisfying the equality $E(f_0,g_0,h_0)(s_0) = 0$. By the continuity of $E$, we can choose a natural number $u$, such that
$E(f,g,h)(s_0) = 0$, whenever $f,g,h\in\mathbb{T}_1$ and $f(t) = f_0(t)$, $g(t) = g_0(t)$,
$h(t) = h_0(t)$ for all $t\le u$. Let $P,Q,R$ be the following  unary operators:
$$P={^{u+1}}f_0,\ Q={^{u+1}}g_0,\ R={^{u+1}}h_0.$$
These operators belong to $\mathbf{O}$ and,
for any $f,g,h\in\mathbb{T}_1$, the
functions $P(f)$, $Q(g)$, $R(h)$ coincide, respectively, with the functions $f_0,g_0,h_0$
on $\{t\in\mathbb{N}\,|\,t\le u\}$ and with the functions $f,g,h$ on $\{t\in\mathbb{N}\,|\,t>u\}$.
We define $U$ as follows:
$$U=\bigcap_{t=0}^u\left\{\xi\in\mathbb{R}\,\big|\,\left|\frac{f_0(t)-g_0(t)}{h_0(t)+1}-\xi\right|<\frac{1}{t+1}\right\}.$$
Then $U$ is a neighbourhood of $\xi_0$, and whenever a triple $(f,g,h)$ names a real
number belonging to $U$, the triple $(P(f),Q(g),R(h))$ also names this number. Now let us define
\begin{align*}
F_0(f,g,h)&=F(P(f),Q(g),R(h),\check{s_0}),\\
G_0(f,g,h)&=G(P(f),Q(g),R(h),\check{s_0}),\\
H_0(f,g,h)&=H(P(f),Q(g),R(h),\check{s_0}).
\end{align*}
By conditions 1 and 3 of Definition \ref{D1}, Corollary \ref{L7} and the fact that $\check{s_0}$ is $\mathbf{O}$-representable, $F_0$, $G_0$, $H_0$ all belong to $\mathbf{O}$.
Let $\xi\in D\cap U$ and $(f,g,h)$ name $\xi$.
Then $(P(f),Q(g),R(h))$ also names $\xi$ and
moreover, $$E(P(f),Q(g),R(h))(s_0) = 0.$$ It follows that the triple
$(F_0(f,g,h),G_0(f,g,h),H_0(f,g,h))$ names $\theta(\xi)$. By Definition \ref{D4},
we obtain that the restriction of $\theta$ to $D\cap U$ is uniformly $\mathbf{O}${\hyp}computable.$\ \ \Box$

\section{Uniform $\mathbf{O}${\hyp}computability of the locally uniformly $\mathbf{O}${\hyp}computable
functions \\with compact domains}\label{S5}

For any $K\in\mathbb{N}$, let $\delta_K$ be the function from $\mathbb{T}_{2K+1}$ defined as follows:\footnote{We will use the functions $\delta_K$ to organize a simple bounded search.}
for all $x_1,y_1,\ldots,$ $x_K,y_K,z$ in $\mathbb{N}$, if $x_i=0$ for some $i\in\{1,\ldots,K\}$ then $$\delta_K(x_1,y_1,\ldots,x_K,y_K,z)=y_i$$ with the least such $i$, otherwise $$\delta_K(x_1,y_1,\ldots,x_K,y_K,z)=z.$$ 
In particular,
$$\delta_1(x,y,z)=\left\{\begin{array}{ll}\!y & \mathrm{if}\ x=0,\\\!z & \mathrm{otherwise.}\end{array}\right.$$

\begin{definition}\label{D7}
A class $\mathbf{O}$ of operators will be called {\em decent} if $\mathbf{O}$ is appropriate and the functions $\lambda x.x+1$,
$\lambda xy.x\dotminus y$ and $\delta_1$ are $\mathbf{O}$-representable.
\end{definition}

\begin{lemma}\label{Orf}
Let $\mathbf{O}$ be a decent class of operators. Then all constant functions from $\mathbb{T}_1$ and the functions $\delta_K$, $K=0,1,2,3,\ldots,$ are $\mathbf{O}$-representable.
\end{lemma}

{\bf Proof.} By using the $\mathbf{O}${\hyp}representability of $\lambda x.x+1$ and $\lambda xy.x\dotminus y$ together with Lemmas \ref{L4} and \ref{L5} we easily see that all constant functions from $\mathbb{T}_1$ are $\mathbf{O}$-representable. The $\mathbf{O}$-representability of the functions $\delta_K$ follows from the equalities
\begin{align*}
\delta_0(z)&=z,\\
\delta_{K+1}(x_1,y_1,x_2,y_2\ldots,x_{K+1},y_{K+1},z)&=\delta_1(x_1,y_1,\delta_K(x_2,y_2,\ldots,x_{K+1},y_{K+1},z))
\end{align*}
by induction on $K$ with application of Lemmas \ref{L4} and \ref{L5}.$\ \ \Box$

\begin{corollary}\label{nc}
If $\mathbf{O}$ is a decent class of continuous operators then all conditionally $\mathbf{O}${\hyp}computable real functions are locally uniformly $\mathbf{O}${\hyp}computable.
\end{corollary}

{\bf Proof.} If $\mathbf{O}$ is a decent class of operators then, by the above lemma, the equality
$$\mu_{k,c}(x,y)=\delta_1(x\dotminus k,\delta_1(k\dotminus x,c,y),y)$$
and Lemmas \ref{L4}, \ref{L5}, the functions listed in the premise of Theorem \ref{T2} are $\mathbf{O}$-representable.$\ \ \Box$

\vskip3mm
The next theorem generalizes and somewhat strengthens Theorem 3 of \cite{s-g}.

\begin{theorem}\label{T3}
Let $\mathbf{O}$ be a decent class of operators. Then all locally uniformly $\mathbf{O}${\hyp}computable real
functions with compact domains are uniformly $\mathbf{O}${\hyp}computable.\footnote{The conclusion of theorem 3 in \cite{s-g} is equivalent to the particular instance of the present statement for the case when $\mathbf{O}$ is the class of the $\cal F${\hyp}substitutional operators for a class $\cal F$ satisfying the conditions of the theorem in question. The class $\mathbf{O}$ in that case will be surely decent by Lemmas \ref{L1}, \ref{L3} and the equality $\,\delta_1(x,y,z)=y(1\dotminus x)+z(1\dotminus (1\dotminus x))$.}
\end{theorem}

{\em Proof}. Suppose $N\in\mathbb{N}$, $\theta:D\to\mathbb{R}$, where $D$ is a compact subset of $\mathbb{R}^N$,
and $\theta$ is locally uniformly $\mathbf{O}${\hyp}computable. Then there exist $K\in\mathbb{N}$, rational
numbers $a_{ij}$ ($i = 1,\ldots,K,\ j = 1,\ldots,N$) and positive rational numbers 
$d_1,\ldots,d_K$, such that $D\subseteq U_1\cup\ldots\cup U_K$, where, for $i = 1,\ldots,K$,
$$U_i = \left\{\,(\xi_1,\ldots,\xi_N)\in\mathbb{R}^N\,\big|\ |\xi_1 - a_{i1}| < d_i,\,\ldots,\,|\xi_N - a_{iN}| < d_i\,\right\}$$
and the restriction of $\theta$ to $D\cap U_i$ is uniformly $\mathbf{O}${\hyp}computable. We will prove
that $\theta$ is also uniformly $\mathbf{O}${\hyp}computable. (Of course, the case $K<2$ is trivial, so we may assume that $K\ge 2$.) In order to prove the uniform $\mathbf{O}${\hyp}computability of~$\theta$, we consider the
continuous function
$$\rho(\xi_1,\ldots,\xi_N) = \max_{i =1,\ldots,K} (d_i - \!\max_{j =1,\ldots,N}\,|\xi_j - a_{ij}|).$$
Since $\rho(\bar{\xi}) > 0$ for all $\bar{\xi}\in D$, there exists a natural number $k$, such that \mbox{$\rho(\bar{\xi})\ge\frac{2}{k+1}$}
for any $\bar{\xi}\in D$. For such a $k$, as it is easy to see, whenever
$(\xi_1,\ldots,\xi_N)\in D$ and $x_1,\ldots,x_N$ are rational numbers satisfying the inequalities $|x_j- \xi_j| < \frac{1}{k+1}$ $(j=1,\ldots,N)$, at least one of the numbers
$$r_1 =d_1 - \!\max_{j =1,\ldots,N}\,|x_j - a_{1j}|,\ \ldots,\ r_K =d_K -\!\max_{j =1,\ldots,N}\, |x_j - a_{Kj}|$$ will be
greater than $\frac{1}{k+1}$,
and $(\xi_1,\ldots,\xi_N)$ will belong to $U_i$ for any $i\in\{1,\ldots,K\}$,
such that $r_i\!>\!\!\frac{1}{k+1}$.
In particular, that will be the case, whenever
\mbox{$(\xi_1,\ldots,\xi_N)\in D$,} $(f_1,g_1,h_1),\,\ldots,\,(f_N,g_N,h_N)$ are triples naming $\xi_1,\,\ldots,\,\xi_N$,
respectively, and
$$x_j=\frac{f_j(k)-g_j(k)}{h_j(k)+1},\ \ j=1,\ldots,N.$$
For all $i = 1,\ldots,K$, let us choose operators $F_i,G_i,H_i\in\mathbf{O}$,
according to Definition \ref{D4}, applied for the restriction of $\theta$ to $D\cap U_i$ 
(which is uniformly $\mathbf{O}${\hyp}computable). We define $3N$-ary operators $F,G,H$ by
\begin{align*}
F(f_1,g_1,h_1,\ldots,f_N,g_N,h_N)&=F_l(f_1,g_1,h_1,\ldots,f_N,g_N,h_N),\\
G(f_1,g_1,h_1,\ldots,f_N,g_N,h_N)&=G_l(f_1,g_1,h_1,\ldots,f_N,g_N,h_N),\\
H(f_1,g_1,h_1,\ldots,f_N,g_N,h_N)&=H_l(f_1,g_1,h_1,\ldots,f_N,g_N,h_N),
\end{align*}
where $l$ is the least of the numbers $i\in\{1,\ldots,K\}$ satisfying the inequality
\begin{equation}\label{E6}
d_i-\!\max_{j =1,\ldots,N}\left|\frac{f_j(k)-g_j(k)}{h_j(k)+1}-a_{ij}\right|>\frac{1}{k+1},
\end{equation}
if there exists such an $i$, and
$$F(f_1,g_1,h_1,\ldots,f_N,g_N,h_N)=G(f_1,g_1,h_1,\ldots,f_N,g_N,h_N)=H(f_1,g_1,h_1,\ldots,f_N,g_N,h_N)=\check{0},$$ 
otherwise.

The above reasoning will show that $F,G,H$ are witnesses for the
uniform $\mathbf{O}${\hyp}computability of $\theta$, if we succeed to prove that they belong to $\mathbf{O}$. Of course, the inequality (\ref{E6}) is equivalent to
\begin{equation}\label{nh}
\max_{j =1,\ldots,N}\left|\frac{f_j(k)-g_j(k)}{h_j(k)+1}-a_{ij}\right|<d_i-\frac{1}{k+1}.
\end{equation}
We will prove the following auxiliary statement.

\vskip2mm
{\em For any rational numbers $a_1,\ldots,a_N,q$, there exists an $\mathbf{O}${\hyp}representable function $e\in\mathbb{T}_{3N}$ such that, for all $x_1,y_1,z_1,\ldots,x_N,y_N,z_N\in\mathbb{N}$, the equivalence
$$e(x_1,y_1,z_1,\ldots,x_N,y_N,z_N)=0\Leftrightarrow\max_{j =1,\ldots,N}\left|\frac{x_j-y_j}{z_j+1}-a_j\right|<q$$
holds.}

\vskip2mm
Then, for any $i\in\{1,\ldots,K\}$, the inequality (\ref{nh}) will be equivalent to some equality of the form
$$e_i(f_1(k),g_1(k),h_1(k),\ldots,f_N(k),g_N(k),h_N(k))=0$$
with $\mathbf{O}${\hyp}representable $e_i$. By the definition of the operator $F$, this will yield the equality
\begin{multline*}
F(f_1,g_1,h_1,\ldots,f_N,g_N,h_N)(t)\\=\delta_K(S_1(f_1,g_1,h_1,\ldots,f_N,g_N,h_N)(t),F_1(f_1,g_1,h_1,\ldots,f_N,g_N,h_N)(t),\ldots,\\S_K(f_1,g_1,h_1,\ldots,f_N,g_N,h_N)(t),F_K(f_1,g_1,h_1,\ldots,f_N,g_N,h_N)(t),0),
\end{multline*}
where
$$S_i(f_1,g_1,h_1,\ldots,f_N,g_N,h_N)(t)=\\ 
e_i(f_1(k),g_1(k),h_1(k),\ldots,f_N(k),\,g_N(k),h_N(k)),\ \ i=1,\ldots,K.$$
Using Corollary \ref{L7} and statements 2, 3 in Lemma \ref{is}, first by the $\mathbf{O}${\hyp}representability of $\check{k}$ and $e_1,\ldots,e_K$ we will be able to conclude that
$S_1,\ldots,S_K \in \mathbf{O}$ and then by the $\mathbf{O}${\hyp}representability of $\delta_K$ and $\check{0}$ and the fact that $F_1,\ldots,F_K \in \mathbf{O}$, it will be true that $F\in\mathbf{O}$. It could be seen in a similar way that $G$ and $H$ also belong to $\mathbf{O}$.

To complete the proof, it remains to prove the auxiliary statement. We will firstly prove that, for any rational number $a$, there exist $\mathbf{O}$-representable ternary functions $\mathrm{lt}_a$ and $\mathrm{gt}_a$ such that, for all $x,y,z\in\mathbb{N}$, the equivalences
$$\mathrm{lt}_a(x,y,z)>0\ \Leftrightarrow\ \frac{x-y}{z+1}<a,\ \ \ \mathrm{gt}_a(x,y,z)>0\ \Leftrightarrow\ \frac{x-y}{z+1}>a$$
hold. Actually, it is sufficient to show how to construct the function
$\mathrm{lt}_a$, since then we may set
$$\mathrm{gt}_a(x,y,z)=\mathrm{lt}_{-a}(y,x,z).$$
If $a=0$ then the inequality
\begin{equation}\label{lt}
\frac{x-y}{z+1}<a
\end{equation}
is equivalent to $x<y$ and we may set
$$\mathrm{lt}_a(x,y,z)=y\dotminus x.$$
To settle the case when $a\ne 0$, we will first construct, for any positive integers $b,c$, an $\mathbf{O}$-representable function $\gamma_{b,c}:\mathbb{N}^{b+c}\to\mathbb{N}$ such that, for all $x_1,\ldots,x_b,y_1,\ldots,y_c\in\mathbb{N}$,
\begin{equation}\label{E9}
\gamma_{b,c}(x_1,\ldots,x_b,y_1,\ldots,y_c)>0\ \Leftrightarrow\ x_1+\cdots+x_b>y_1+\cdots+y_c.
\end{equation}
The construction is by the following inductive definition, where $b$ and $c$ can be arbitrary positive integers:
\begin{align*}
\gamma_{1,1}(x_1,y_1)&=x_1\dotminus y_1,\\ 
\gamma_{1,c+1}(x_1,y_1,\ldots,y_c,y_{c+1})&=\gamma_{1,c}(x_1\dotminus y_{c+1},y_1,\ldots,y_c),\\ 
\gamma_{b+1,c+1}(x_1,\ldots,x_b,x_{b+1},y_1,\ldots,y_c,y_{c+1})&=\delta_1(x_{b+1}\dotminus y_{c+1},\\\gamma_{b,c+1}&(x_1,\ldots,x_b,y_1,\ldots,y_c,y_{c+1}\dotminus x_{b+1}),\\\gamma_{b+1,c}&(x_1,\ldots,x_b,x_{b+1}\dotminus y_{c+1},y_1,\ldots,y_c)).
\end{align*}
One proves inductively that all functions $\gamma_{b,c}$ satisfy the equivalence (\ref{E9}) and are $\mathbf{O}$-representable (of course, the $\mathbf{O}$-representability of $\lambda xy.x\dotminus y$ and $\delta_1$, as well as Lemmas \ref{L4} and \ref{L5} are used in the proof of the last statement). If $a>0$ then $a=b/c$ with some positive integers $b,c$, and (\ref{lt}) is equivalent to the inequality $c.(x-y)<b.(z+1)$. The last inequality is easily seen to be equivalent to $c.(x\dotminus y)<b.(z+1)$. Thus we may set 
$$\mathrm{lt}_a(x,y,z)=\gamma_{b,c}(\underbrace{z+1,\ldots,z+1}_{b\text{ times}},\underbrace{x\dotminus y,\ldots,x\dotminus y}_{c\text{ times}})$$
in this case. Finally, if $a<0$ then $a=-c/b$ with some positive integers $b,c$, and~(\ref{lt}) is equivalent to the inequality $b.(y-x)>c.(z+1)$, which, in turn, is equivalent to $b.(y\dotminus x)>c.(z+1)$. Therefore we may set
$$\mathrm{lt}_a(x,y,z)=\gamma_{b,c}(\underbrace{y\dotminus x,\ldots,y\dotminus x}_{b\text{ times}},\underbrace{z+1,\ldots,z+1}_{c\text{ times}})$$
now. By using the $\mathbf{O}${\hyp}representability of the functions $\lambda x.x+1$, $\lambda xy.x\dotminus y$ and $\gamma_{b,c}$ ($b,c=1,2,3,\ldots$), and making use again of Lemmas \ref{L4} and \ref{L5}, we see that all the functions $\mathrm{lt}_a$ defined above are $\mathbf{O}$-representable. 

Suppose now that rational numbers $a_1,\ldots,a_N,q$ are given. The inequality in the right-hand side of the equivalence in the auxiliary statement is equivalent to the conjunction of the inequalities
$$gt_{b_j}(x_j,y_j,z_j)>0,\ \ lt_{c_j}(x_j,y_j,z_j)>0,\ \ j=1,\ldots,N,$$
where $b_j=a_j-q$, $c_j=a_j+q$. Let us set
$$u_j(x,y,z)=\delta_1(gt_{b_j}(x,y,z),gt_{b_j}(x,y,z),lt_{c_j}(x,y,z)),\ \ j=1,\ldots,N,$$
and then set
$$e(x_1,y_1,z_1,\ldots,x_N,y_N,z_N)=\delta_N(u_1(x_1,y_1,z_1),1,\ldots,u_N(x_N,y_N,z_N),1,0).$$
It is easy to check the equivalence from the auxiliary statement, and the $\mathbf{O}${\hyp}representability of the function $e$ follows from the $\mathbf{O}${\hyp}representability of the functions $gt_{b_j},lt_{c_j},\check{1},\check{0},\delta_1,\delta_N$ and Lemmas \ref{L4}, \ref{L5}.$\ \ \Box$ 

\begin{corollary}
If $\mathbf{O}$ is a decent class of continuous operators, then all conditionally $\mathbf{O}${\hyp}computable real functions with compact
domains are uniformly $\mathbf{O}${\hyp}computable.
\end{corollary}

{\bf Proof.} By Corollary \ref{nc} and the above theorem.$\ \ \Box$

\section{Appendix: Conditional computability of functions in effective metric spaces}
The referee of the paper asked if Definitions \ref{D4}, \ref{D5} and \ref{D6} can be extended to functions between metric spaces other than reals, for instance using the representation-theoretic approach of TTE \cite[Definition 8.1.2]{w}.

\vskip2mm
To do such an extension, we have to consider effective metric spaces $\mathbf{M}=(M,d,A,\alpha)$ and $\mathbf{M^\prime}=(M^\prime,d^\prime,A^\prime,\alpha^\prime)$ in the sense of the above-mentioned definition instead of $\mathbb{R}^N$ and $\mathbb{R}$, respectively. However, we must assume that the domains of $\alpha$ and $\alpha^\prime$ consist of natural numbers rather than strings over an arbitrary finite alphabet -- this is needed, since our operators act on functions in $\mathbb{N}$. An {\em ordinary name} of an element $\xi$ of $M$ will be, by definition, a total one-argument function $f$ in $\mathbb{N}$ such that $f(t)\in\mathrm{dom}(\alpha)$ and $d(\alpha(f(t)),\xi)<\frac{1}{t+1}$ for any $t\in\mathbb{N}$ (similarly for ordinary names of the elements of $M^\prime$).\footnote{We add the adjective ``ordinary'' in front of ``name'' in order to distinguish the names used here from the Cauchy names used in \cite{w}.} The next definition contains analogs of the above-mentioned definitions for the case of effective metric spaces of the above sort.

\begin{definition} 
Let $\mathbf{O}$ be a class of operators, let $\mathbf{M}=(M,d,A,\alpha)$, $\mathbf{M^\prime}=(M^\prime,d^\prime,A^\prime,\alpha^\prime)$ be effective metric spaces with $\mathrm{dom}(\alpha),\mathrm{dom}(\alpha^\prime)\subseteq\mathbb{N}$, and let $\theta:D\to M^\prime$, where $D\subseteq M$. The function $\theta$ will be called {\em uniformly computable} with respect to~$\mathbf{O}$
or {\em uniformly $\mathbf{O}${\hyp}computable}, for short, if there exists a unary operator $T\in\mathbf{O}$ such that, whenever $\xi\in D$ and $f$ is an ordinary name of $\xi$, the function $T(f)$ is an ordinary name of $\theta(\xi)$. The function $\theta$ will be called {\em conditionally computable} with respect to~$\mathbf{O}$ or {\em conditionally $\mathbf{O}${\hyp}computable}, for short, if there exist a unary operator~$E\in\mathbf{O}$ and a binary operator $T\in\mathbf{O}$, such that, whenever $\xi\in D$ and $f$ is an ordinary name of $\xi$, there exists a natural number $s$ satisfying the equality $E(f)(s)=0$, and the function $T(f,\check{s})$ is an ordinary name of $\theta(\xi)$ for any such $s$. The function $\theta$ will be called {\em locally uniformly $\mathbf{O}${\hyp}computable}, if any point of $D$ has some neighbourhood $U$ in $\mathbf{M}$ such that the restriction of $\theta$ to $D\cap U$ is uniformly $\mathbf{O}${\hyp}computable.
\end{definition}

Here are some analogs of Theorems \ref{T1}, \ref{T2} and \ref{T3}.

\begin{theorem}[Analog of Theorem \ref{T1}]
Let $\mathbf{O}$ be an appropriate class of operators. Let a two-argument $\mathbf{O}$-representable function $C$ and $\mathbf{O}${\hyp}representable one-argument functions $L$ and $R$ in~$\mathbb{N}$ exist such that 
$$\{(u,v)\in\mathbb{N}^2\,|\,C(u,v)=0\}=\{(0,0)\},\ \ \{(L(s),R(s))\,|\,s\in\mathbb{N}\} = \mathbb{N}^2.$$ Let $\mathbf{M_0}=(M_0,d_0,A_0,\alpha_0),\ \mathbf{M_1}=(M_1,d_1,A_1,\alpha_1),\ \mathbf{M_2}=(M_2,d_2,A_2,\alpha_2)$ be effective metric spaces with $\mathrm{dom}(\alpha_0),\mathrm{dom}(\alpha_1),\mathrm{dom}(\alpha_2)\subseteq\mathbb{N}$, $\theta_0$ be a conditionally $\mathbf{O}${\hyp}computable partial function from $\mathbf{M_1}$ to $\mathbf{M_0}$, and $\theta_1$ be a conditionally $\mathbf{O}${\hyp}computable partial function from $\mathbf{M_2}$ to $\mathbf{M_1}$. Then the composition of $\theta_0$ and $\theta_1$ is also conditionally
$\mathbf{O}${\hyp}computable.
\end{theorem}

{\bf Proof.} Let $\theta$ be the composition of $\theta_0$ and $\theta_1$, i.e. $\theta$ is the partial function from $\mathbf{M_2}$ to $\mathbf{M_0}$ defined by $\theta(\xi) = \theta_0(\theta_1(\xi))$. For $i = 0,1$,
let $E_i,T_i$ be operators from $\mathbf{O}$, such that $$\exists s(E_i(f)(s) = 0)\ \&\
\forall s((E_i(f)(s) = 0)
\Rightarrow T_i(f,\check{s})\
\textrm{is an ordinary name of}\ \theta_i(\xi))$$
for any $\xi\in\mathrm{dom}(\theta_i)$ and any ordinary name $f$ of $\xi$. Let the operators $E$ and $T$ be defined as follows:
\begin{align*}
E(f)(s)=&\ C(E_1(f)(R(s)),
E_0(T_1(f,\mathring{R}(\check{s})))(L(s))),\\
T(f,e)=&\ T_0(T_1(f,\mathring{R}(e)),\mathring{L}(e)).
\end{align*}
The reasoning continues in 
the same manner, as in the proof of Theorem \ref{T1}.
It is seen that the operators $E$ and $T$ are witnesses for the conditional $\mathbf{O}${\hyp}computability of $\theta$.$\ \ \Box$

\begin{theorem}[Analog of Theorem \ref{T2}]
Let $\mathbf{O}$ be an appropriate class of continuous operators, and let the functions $\check{c}$ for all $c\in\mathbb{N}$ and the functions
$$\mu_{k,c}(x,y)=\left\{\begin{array}{ll}\!c & \mathrm{if}\ x=k,\\\!y & \mathrm{otherwise}\end{array}\right.$$
for all $k,c\in\mathbb{N}$ be $\mathbf{O}$-representable. Let $\mathbf{M}=(M,d,A,\alpha)$, $\mathbf{M^\prime}=(M^\prime,d^\prime,A^\prime,\alpha^\prime)$ be effective metric spaces with $\mathrm{dom}(\alpha),\mathrm{dom}(\alpha^\prime)\subseteq\mathbb{N}$. Then all conditionally $\mathbf{O}${\hyp}computable partial functions from $\mathbf{M}$ to $\mathbf{M^\prime}$ are locally uniformly $\mathbf{O}${\hyp}computable.
\end{theorem}

{\bf Proof.} Let $\theta:D\to M^\prime$, where $D\subseteq M$, be a conditionally $\mathbf{O}${\hyp}computable partial function,
and let $\xi_0\in D$. Let $E$ and $T$ be operators from $\mathbf{O}$, such that $$\exists s(E(f)(s) = 0)\ \&\
\forall s((E(f)(s) = 0)
\Rightarrow T(f,\check{s})\
\textrm{is an ordinary name of}\ \theta(\xi))$$
for any $\xi\in\mathrm{dom}(\theta)$ and any ordinary name $f$ of $\xi$.
Let $f_0$ be an ordinary name of $\xi_0$, and let $s_0$ be a natural number,
satisfying the equality $E(f_0)(s_0) = 0$. By the continuity of $E$, we can choose a natural number $u$, such that
$E(f)(s_0) = 0$, whenever $f\in\mathbb{T}_1$ and $f(t) = f_0(t)$ for all $t\le u$. Let $P={^{u+1}}f_0$ (in the notations used in the proof of Theorem \ref{T2}), let $U$ be defined by
$$U=\bigcap_{t=0}^u\left\{\xi\in M\,\big|\,d(\alpha(f_0(t)),\xi)<\frac{1}{t+1}\right\}$$
and let us set $$T_0(f)=T(P(f),\check{s_0}).$$
Then, as in the proof of Theorem \ref{T2}, $U$ is a neighbourhood of $\xi_0$, $P$ belongs to $\mathbf{O}$
and the operator $T_0$ is a witness for the uniform $\mathbf{O}${\hyp}computability of the restriction
of $\theta$ to $D\cap U$.$\ \ \Box$

\begin{theorem}[Analog of Theorem \ref{T3}]
Let $\mathbf{O}$ be an appropriate class of operators such that the function $\delta_1$ and all constant functions from $\mathbb{T}_1$ are $\mathbf{O}${\hyp}representable. Let $\mathbf{M}=(M,d,A,\alpha)$, $\mathbf{M^\prime}=(M^\prime,d^\prime,A^\prime,\alpha^\prime)$ be effective metric spaces with $\mathrm{dom}(\alpha),\mathrm{dom}(\alpha^\prime)$ $\subseteq\mathbb{N}$, and let, for any $a$ in $A$ and any rational number $q$, there exists an $\mathbf{O}${\hyp}representable function from $\mathbb{T}_1$ having the value 0 exactly for those $n\in\mathrm{dom}(\alpha)$ which satisfy the inequality $d(\alpha(n),a)<q$. Then all locally uniformly $\mathbf{O}${\hyp}computable partial
functions with compact domains from $\mathbf{M}$ to $\mathbf{M^\prime}$ are uniformly $\mathbf{O}${\hyp}computable.
\end{theorem}

{\em Proof}. Let $\theta:D\to\mathbb{M^\prime}$ be a locally uniformly $\mathbf{O}${\hyp}computable partial 
function from $\mathbf{M}$ to $\mathbf{M^\prime}$, and let its domain $D$ be compact. Then there exist $K\in\mathbb{N}$, elements $a_1,\ldots,a_K$ of $A$ and positive rational numbers $d_1,\ldots,d_K$
such that $D\subseteq U_1\cup\ldots\cup U_K$, where, for $i = 1,\ldots,K$,
$$U_i = \{\,\xi\,|\,d(\xi,a_i)<d_i\}$$
and the restriction of $\theta$ to $D\cap U_i$ is uniformly $\mathbf{O}${\hyp}computable. We will prove
that $\theta$ is also uniformly $\mathbf{O}${\hyp}computable. (Of course, the case $K<2$ is trivial, so we may assume that $K\ge 2$.) In order to prove the uniform $\mathbf{O}${\hyp}computability of~$\theta$, we consider the continuous function
$$\rho(\xi) = \max_{i =1,\ldots,K}\,(d_i - d(\xi,a_i)).$$
Since $\rho(\xi) > 0$ for all $\xi\in D$, there exists a natural number $k$, such that \mbox{$\rho(\xi)\ge\frac{2}{k+1}$}
for any $\xi\in D$ and let us choose such a $k$. For any $\xi \in D$ and any ordinary name $f$ of $\xi$,
as in the proof of Theorem \ref{T3}, at least one of the numbers
$$r_1 =d_1 - d(\alpha(f(k)),a_1),\ \ldots,\ r_K =d_K - d(\alpha(f(k)),a_K)$$ 
will be greater than $\frac{1}{k+1}$
and $\xi$ will belong to $U_i$ for any $i\in\{1,\ldots,K\}$,
such that $r_i\!>\!\!\frac{1}{k+1}$.
For any $i\in\{1,\ldots,K\}$, let us choose an operator $T_i\in\mathbf{O}$, such that, whenever $\xi\in D\cap U_i$ and $f$ is an ordinary name of $\xi$, the function $T_i(f)$ is an ordinary name of $\theta(\xi)$
(such operators exist due to the uniform $\mathbf{O}${\hyp}computability of the restriction of $\theta$ to any of the sets $D\cap U_1,\ldots,D\cap U_K$). For any such $i$, let $e_i$ be an $\mathbf{O}${\hyp}representable function from~$\mathbb{T}_1$ such that
$$\forall n\in\mathrm{dom}(\alpha)\left(e_i(n)=0\Leftrightarrow d(\alpha(n),a_i)<d_i-\frac{1}{k+1}\right).$$
We define an operator $T$ by setting
$$T(f)=T_l(f),$$
where $l$ is the least of the numbers $i\in\{1,\ldots,K\}$ satisfying the equality $e_i(f(k))=0$ if there exists such an $i$, and
$$T(f)=\check{0},$$
otherwise. The above reasoning will show that $T$ is a witness for the
uniform $\mathbf{O}${\hyp}computability of $\theta$ if we succeed to prove that $T\in\mathbf{O}$.
To show this, we note that, for any $f\in\mathbb{T}_1$, the equality
$$T(f)(t)=\delta_K(e_1(f(k)),T_1(f)(t),\ldots,e_K(f(k)),T_K(f)(t),0)$$
holds. As in the proof of Lemma \ref{Orf}, we can show that all functions $\delta_1,\delta_2,\delta_3,\ldots$ are $\mathbf{O}${\hyp}representable. The $\mathbf{O}${\hyp}representability of the functions $\delta_K,e_1,\ldots,e_K,\check{k}$, together with Corollary~\ref{L7} and statements 2, 3 in Lemma \ref{is}, imply that $T\in\mathbf{O}.\ \ \Box$ 

Unfortunately, we need an additional assumption about the class $\mathbf{O}$ for being able to consider Theorems \ref{T1}, \ref{T2} and \ref{T3} as particular instances of their analogs proved above. The assumption is, roughly speaking, about the existence of an $\mathbf{O}${\hyp}representable mechanism in $\mathbb{N}$ for coding of ordered pairs and the corresponding decoding. Namely, let $\mathbf{O}$ be an appropriate class of operators, and let there exist an $\mathbf{O}${\hyp}representable function $J\in\mathbb{T}_2$ and $\mathbf{O}${\hyp}representable functions $L,R\in\mathbb{T}_1$ such that $L(J(u,v))=u$ and $R(J(u,v))=v$ for all $u,v\in\mathbb{N}$. Then we can turn any $\mathbb{R}^N$ into an effective metric spaces $\mathbb{M}_N$ in such a way that the notions introduced Definitions \ref{D4}, \ref{D5} and \ref{D6} are equivalent to the corresponding notions introduced here for partial functions from $\mathbb{M}_N$ to $\mathbb{M}_1$.

First of all, the above assumptions imply that, for any positive integer $K$ there exist an $\mathbf{O}${\hyp}representable function $J_K\in\mathbb{T}_K$ and $\mathbf{O}${\hyp}representable functions $P_{K,1},\ldots,P_{K,K}\in\mathbb{T}_1$ such that
$$P_{K,i}(J_K(u_1,\ldots,u_K))=u_i,\ \ i=1,\ldots,K,$$
for all $u_1,\ldots,u_K\in\mathbb{N}$.

For any positive integer $N$, we consider the effective metric space $\mathbb{M}_N=(\mathbb{R}^N,d_N,\mathbb{Q}^N,\alpha_N)$, where, as usually, $\mathbb{Q}$ is the set of the rational numbers, and $d_N,\alpha_N$ are defined as follows:
\begin{gather*}
d_N((\xi_1,\ldots,\xi_N),(\eta_1,\ldots,\eta_N))=\max(|\xi_1-\eta_1|,\ldots,|\xi_N-\eta_N|),\\
\alpha_N(n)=\left(\frac{P_{3N,1}(n)-P_{3N,2}(n)}{P_{3N,3}(n)+1},\ldots,\frac{P_{3N,3N-2}(n)-P_{3N,3N-1}(n)}{P_{3N,3N}(n)+1}\right).
\end{gather*}
It is easy to prove that a partial function from $\mathbb{R}^N$ to $\mathbb{R}$ is uniformly or conditionally $\mathbf{O}${\hyp}computable if and only if it is uniformly or conditionally $\mathbf{O}${\hyp}computable, respectively, as a partial function from $\mathbb{M}_N$ to $\mathbb{M}_1$ (after the identification of $\mathbb{R}^1$ and $\mathbb{Q}^1$ with $\mathbb{R}$ and $\mathbb{Q}$, respectively). Namely, if $\theta:D\to\mathbb{R}$, where $D\subseteq\mathbb{R}^N$, then:
\begin{enumerate}
\item If $F,G,H$ have the properties from Definition \ref{D4}, then the operator $T$ defined by
\begin{align*}
T(f)=\mathring{J_3}(F(\mathring{P_{3N,1}}(f),\mathring{P_{3N,2}}(f)&,\mathring{P_{3N,3}}(f),\ldots,\mathring{P_{3N,3N-2}}(f),\mathring{P_{3N,3N-1}}(f),\mathring{P_{3N,3N}}(f)),\\G(\mathring{P_{3N,1}}(f),\mathring{P_{3N,2}}(f)&,\mathring{P_{3N,3}}(f),\ldots,\mathring{P_{3N,3N-2}}(f),\mathring{P_{3N,3N-1}}(f),\mathring{P_{3N,3N}}(f)),\\H(\mathring{P_{3N,1}}(f),\mathring{P_{3N,2}}(f)&,\mathring{P_{3N,3}}(f),\ldots,\mathring{P_{3N,3N-2}}(f),\mathring{P_{3N,3N-1}}(f),\mathring{P_{3N,3N}}(f))) 
\end{align*}
is a witness for the $\mathbf{O}${\hyp}computability of $\theta$ as a partial function from $\mathbb{M}_N$ to~$\mathbb{M}_1$.
\item If $T$ is a witness for the $\mathbf{O}${\hyp}computability of $\theta$ as a partial function from $\mathbb{M}_N$ to~$\mathbb{M}_1$, then the operators $F,G,H$ defined by
\begin{align*}
F(f_1,g_1,h_1,\ldots,f_N,g_N,h_N)\!=\,&\mathring{P_{3,1}}(T(\mathring{J_{3N}}(f_1,g_1,h_1,\ldots,f_N,g_N,h_N))),\\
G(f_1,g_1,h_1,\ldots,f_N,g_N,h_N)\!=\,&\mathring{P_{3,2}}(T(\mathring{J_{3N}}(f_1,g_1,h_1,\ldots,f_N,g_N,h_N))),\\
H(f_1,g_1,h_1,\ldots,f_N,g_N,h_N)\!=\,&\mathring{P_{3,3}}(T(\mathring{J_{3N}}(f_1,g_1,h_1,\ldots,f_N,g_N,h_N)))
\end{align*}
have the properties from Definition \ref{D4}.
\item If the operators $E,F,G,H$ have the properties from Definition \ref{D5}, then the operators $E^\prime$ and $T$ defined by
\begin{align*}
E^\prime(f)=E(\mathring{P_{3N,1}}(f),\mathring{P_{3N,2}}(f)&,\mathring{P_{3N,3}}(f),\ldots,\mathring{P_{3N,3N-2}}(f),\mathring{P_{3N,3N-1}}(f),\mathring{P_{3N,3N}}(f)),\\ 
T(f,e)=\mathring{J_3}(F(\mathring{P_{3N,1}}(f),\mathring{P_{3N,2}}&(f),\mathring{P_{3N,3}}(f),\ldots,\mathring{P_{3N,3N-2}}(f),\mathring{P_{3N,3N-1}}(f),\mathring{P_{3N,3N}}(f),e),\\G(\mathring{P_{3N,1}}(f),\mathring{P_{3N,2}}&(f),\mathring{P_{3N,3}}(f),\ldots,\mathring{P_{3N,3N-2}}(f),\mathring{P_{3N,3N-1}}(f),\mathring{P_{3N,3N}}(f),e),\\H(\mathring{P_{3N,1}}(f),\mathring{P_{3N,2}}&(f),\mathring{P_{3N,3}}(f),\ldots,\mathring{P_{3N,3N-2}}(f),\mathring{P_{3N,3N-1}}(f),\mathring{P_{3N,3N}}(f),e)) 
\end{align*}
are witnesses for the conditional $\mathbf{O}${\hyp}computability of $\theta$ as a partial function from $\mathbb{M}_N$ to~$\mathbb{M}_1$.
\item If the operators $E^\prime$ and $T$ are witnesses for the conditional $\mathbf{O}${\hyp}computability of $\theta$ as a partial function from $\mathbb{M}_N$ to~$\mathbb{M}_1$, then the operators $E,F,G,H$ defined by
\begin{align*}
E(f_1,g_1,h_1,\ldots,f_N,g_N,h_N)\!=\,&E^\prime(\mathring{J_{3N}}(f_1,g_1,h_1,\ldots,f_N,g_N,h_N)),\\
F(f_1,g_1,h_1,\ldots,f_N,g_N,h_N,e)\!=\,&\mathring{P_{3,1}}(T(\mathring{J_{3N}}(f_1,g_1,h_1,\ldots,f_N,g_N,h_N),e)),\\
G(f_1,g_1,h_1,\ldots,f_N,g_N,h_N,e)\!=\,&\mathring{P_{3,2}}(T(\mathring{J_{3N}}(f_1,g_1,h_1,\ldots,f_N,g_N,h_N),e)),\\
H(f_1,g_1,h_1,\ldots,f_N,g_N,h_N,e)\!=\,&\mathring{P_{3,3}}(T(\mathring{J_{3N}}(f_1,g_1,h_1,\ldots,f_N,g_N,h_N),e))
\end{align*}
have the properties from Definition \ref{D5}.
\end{enumerate}

For the derivation of Theorem \ref{T1} from its analog in the case of substitution in a function of more than one argument, the following statement can be additionally proved and used.

\vskip2mm
{\em Let an $\mathbf{O}${\hyp}representable function $C$ with the property from Theorem \ref{T1} exist, and let $\theta_1,\ldots,\theta_K$ be conditionally $\mathbf{O}${\hyp}computable partial functions from an effective metric space $\mathbf{M}=(M,d,A,\alpha)$ to $\mathbb{M}_1$, where $\mathrm{dom}(\alpha)\subseteq\mathbb{N}$. Let the partial function $\theta$ from $\mathbf{M}$ to $\mathbb{M}_K$ be defined by
$$\theta(\xi)=(\theta_1(\xi),\ldots,\theta_K(\xi)).$$
Then $\theta$ is also conditionally $\mathbf{O}${\hyp}computable.}

\vskip2mm
For the derivation of Theorem \ref{T3} from its analog, the following statement should be additionally proved.

\vskip2mm
{\em Let the class $\mathbf{O}$ be decent. Then, for any positive integer $N$, any $a$ in $\mathbb{Q}^N$ and any rational number $q$, there exists an $\mathbf{O}${\hyp}representable function from $\mathbb{T}_1$ having the value 0 exactly for those $n\in\mathbb{N}$ which satisfy the inequality $d_N(\alpha_N(n),a)<q$.}

\vskip2mm
This statement easily follows from the auxiliary statement in the proof of Theorem \ref{T3}.

\section*{Acknowledgments}
We thank the anonymous referee for many useful suggestions.

\bibliographystyle{elsarticle-num}

\end{document}